\documentclass[11pt]{article}
\usepackage{latexsym}
\usepackage{amsfonts}

\newcommand{\ba}{\begin{array}}\newcommand{\ea}{\end{array}}
\newcommand{\sym}{S}
\newcommand{\ns}{\rm}

\newcommand{\nse}{\kern-3pt\ns$=$}\newcommand{\qd}{\hfill$\Box$\medbreak}

\newcommand{\ext}{\raise1pt\hbox{$\ts\bigwedge$}}
\renewcommand{\sym}{S}
\newcommand{\ts}{\textstyle}
\newcommand{\rf}[1]{(\ref{#1})}
\newcommand{\chii}{\raise2pt\hbox{$\chi$}}

\newcommand{\Fg}{\mbox{${\cal F}\kern-2pt_g$}}

\newcommand{\Mg}{\mbox{${\cal M}\kern-2pt_g$}}
\newcommand{\Ng}{\mbox{${\cal N}\kern-2pt_g$}}
\newcommand{\V}{V\kern-1pt}

\newcommand{\Gg}{\mbox{${\cal G}\kern-2pt_g$}}
\newcommand{\cir}{\raise1.6pt\hbox{\footnotesize$\circ$}}
\newcommand{\tr}{\mbox{\bf tr}}

\newcommand{\Res}[2]{\hbox{\ns
Res}\kern-16pt\lower5pt\hbox{\footnotesize$_{#1}$}\kern2pt\left[#2\right]}
\newcommand{\qk}{quaternion-K\"ahler\kern2pt}\renewcommand{\,}{\kern1pt}

\newcommand{\dirac}{/\kern-5pt\partial}

\renewcommand{\ts}{\textstyle}
0
\newtheorem{theo}{Theorem}[section]

\newtheorem{corol}{Corollary}[section]
\def\frac#1#2{{#1\over#2}}

\def\be#1\ee{\begin{equation}#1\end{equation}}

\begin{document}
\title{A formula for the trace of symmetric powers of matrices}

\author{Jose Luis Cisneros\footnote{Instituto de Matem\'aticas, UNAM, Unidad Cuernavaca, A.P.
6--60, C.P. 62131, Cuernavaca, Morelos, M\'exico. E-mail: jlcm@matcuer.unam.mx}, 
Rafael Herrera\footnote{Centro de
Investigaci\'on en Matem\'aticas, A. P. 402,
Guanajuato, Gto., C.P. 36000, M\'exico. E-mail: rherrera@cimat.mx}
\footnote{Partially supported by 
grants of CONACyT and LAISLA (CONACyT-CNRS)} \,\,
and Noemi Santana\footnote{Instituto de Matem\'aticas, UNAM, Unidad Cuernavaca, A.P.
6--60, C.P. 62131, Cuernavaca, Morelos, M\'exico. E-mail: noemi.santana@im.unam.mx}
\footnote{Partially
supported by 
grants of CONACyT and LAISLA (CONACyT-CNRS).}
 }

\date{{\ns \today}}

\maketitle

\vspace{-20pt}

{
\abstract{

We present a formula for the trace of any symmetric power of a $n\times n$ matrix (with
coefficients in a field) in terms of the ordinary powers of the matrix, 
an arbitrarily chosen linear function which
vanishes on the identity matrix, and $n-2$ polynomial functions defined recursively.

}
}

\section{Introduction}

This paper answers a question posed in \cite{Cisneros} regarding the trace of symmetric powers of
a square matrix of arbitrary size.
More precisely, in \cite{Cisneros} it was proved that given a matrix 
$A=\left( \begin{array}{cc}
a & b\\
c & d
\end{array}\right)$
with $b\not=0$, and it $k$-th power
$A^k=\left( \begin{array}{cc}
a_k & b_k\\
c_k & d_k
\end{array}\right)$,
the quantity $b_{k+1}/b$ is invariant under conjugation, and furthermore, 
\[{b_{k+1}\over b} =  \tr (\sym^kA),\]
where $\sym^kA$ is the $k$-th symmetric power of $A$, that is, the matrix giving the induced
transformation on homogeneous polynomials of degree $k$.

Thus, the natural question is the following: can this result be generalized to $n\times n$
matrices by using entries of powers of $A$? Here, we answer this question in the affirmative in a slightly
modified form (see
Theorem \ref{thm:general-case}). The main of the theorem is that the traces of arbitrary symetric powers
of arbitrary matrices can be calculated by using a finite number of polynomial functions.

The paper is organized as follows.
In Section \ref{sec:preliminaries}, we recall the two identities we need, the Cayley-Hamilton
theorem and a $K$-theory formula, and also spot their similarity. In Section \ref{sec:proof}, we define
the functions we need to state Theorem \ref{thm:general-case}, and prove it.
We also give some corollaries. In Section \ref{sec:examples} we give the explicit expressions for
the various functions involved in the formulas for the cases $n=3,4$. 

\vspace{.1in}

{\em Acknowledgements}. The second named author would like to thank the International
Centre for Theoretical Physics (Italy) and the Institut des Hautes \'Etudes Scientifiques (France)
for their hospitality and support.

\section{Preliminaries}\label{sec:preliminaries}

Let $V$ be a vector space over a field $K$,
$\sym^k V$ and $\ext^k V$ denote its $k$-th symmetric and antisymmetric tensor powers respectively.

\subsection{Cayley-Hamilton theorem}
The characteristic polynomial of a $n\times n$ matrix $A$ can be
expressed as
\[p_A(t)=\sum_{i=0}^n (-1)^i\tr(\ext^i A)\, t^{n-i},\]
where $\tr(\ext^i A)$ is the trace of the $i$-th exterior power of $A$.
The Cayley-Hamilton theorem states that the 
matrix $A$ satisfies its own characteristic polynomial, i.e.
\[0=\sum_{i=0}^n (-1)^i\tr(\ext^i A)\,A^{n-i}.\]
This matrix identity can be multiplied by any power of $A$
\begin{equation}
0=\sum_{i=0}^n (-1)^i\tr(\ext^i A)\,A^{n-i+p}.\label{eq:Cayley} 
\end{equation}

\subsection{K-theory formula}
On the other hand, consider the K-theory formula \cite[p.119]{Atiyah}
\[\ext^k(V-W) = \sum_{i=0}^k (-1)^i \ext^{k-i}V\otimes \sym^iW.\]
We will apply this formula when $W=V$. Recall that
$\ext^{n+i} V =0$ for $i\geq 1$ so that
\begin{eqnarray*}
0
&=& \ext^{n+p} (V-V)\\
&=&\sum_{i=0}^{n+p} (-1)^i \ext^{n+p-i}V\otimes \sym^iV\\ 
&=&\ext^{n+p}V - \ext^{n+p-1}V\otimes V + \ext^{n+p-2}V\otimes \sym^2V +\ldots 
+ (-1)^{n+p-1} V\otimes \sym^{n+p-1}V + (-1)^{n+p}  \sym^{n+p}V \\
&=&(-1)^p\ext^{n}V\otimes\sym^pV + (-1)^{p+1}\ext^{n-1}V\otimes \sym^{p+1}V 
+\ldots 
+ (-1)^{n+p-1} V\otimes \sym^{n+p-1}V + (-1)^{n+p}  \sym^{n+p}V \\
&=& (-1)^p\sum_{i=0}^n (-1)^i\ext^{n-i}V\otimes \sym^{p+i}V
\end{eqnarray*}
or equivalently
\begin{eqnarray*}
0
&=&\sum_{j=0}^n (-1)^j \sym^{n+p-j}V \otimes \ext^{j}V.
\end{eqnarray*}
This identity goes through to the induced operators $\sym^k A$ and $\ext^k A$ on the 
symmetric and antisymmetric powers $\sym^k V$ and $\ext^k V$, 
and is preserved by the trace
\begin{equation}
0=\sum_{j=0}^n (-1)^j \tr(\sym^{n+p-j}A)  \tr(\ext^{j}A).\label{eq:K-theory-trace}
\end{equation}

\subsection{Similarity}
Here we notice the similarity between the two identities \rf{eq:Cayley}
\[0=\sum_{j=0}^n (-1)^j\tr(\ext^j A)\,A^{n+p-j}\]
and \rf{eq:K-theory-trace}
\begin{eqnarray*}
0
&=&\sum_{j=0}^n (-1)^j  \tr(\ext^{j}A) \,\,\tr(\sym^{n+p-j}A),
\end{eqnarray*}
whose coefficients are
\[(-1)^j  \tr(\ext^{j}A).\]
Thus, we will use these two formulas to prove the main theorem.

\section{The formulae}\label{sec:proof}
Let $A$ be an arbitrary $n\times n$ matrix with coefficients in a field $K$, $A\in M_{n\times n}(K)$.
Let $\sigma_1$ be any linear functional on $M_{n\times n}(K)$ which vanishes on the identity matrix.
We define $n-2$ polynomial functions $\sigma_i$, $i=2,\ldots, n-1$, in the following recursive way:
\begin{equation}
\left\{
\begin{array}{ccl}
\sigma_1(A^2) & = & \sigma_1(A)\tr(A)-\sigma_2(A)\\
\sigma_1(A^3) & = & \sigma_1(A)\tr(\sym^{2} A)-\sigma_2(A)\tr(A)+\sigma_3 (A)\\
\vdots &  & \hspace{1in}\vdots\\
\sigma_1(A^{n-1}) & = & \sigma_1(A)\tr(\sym^{n-2} A)-\sigma_2(A)\tr(\sym^{n-3}A)+\cdots+(-1)^{n-2}\sigma_{n-1}
(A) \\
\end{array}
\right.\label{eq:sigmas}
\end{equation}

Now, we have $n-1$ identities 
{\footnotesize
\[
\left\{
\begin{array}{rcc}
\tr(\sym^{k} A) - \tr(A)\tr(\sym^{k-1} A) + \tr(\ext^2 A)\tr(\sym^{k-2} A) +\cdots   
+(-1)^n\tr(\ext^n A)\tr(\sym^{k-n} A) & = & 0\\
\tr(\sym^{k-1} A) - \tr(A)\tr(\sym^{k-2} A) + \tr(\ext^2 A)\tr(\sym^{k-3} A) +\cdots  
+(-1)^n\tr(\ext^n A)\tr(\sym^{k-n-1} A) & = & 0\\
\vdots\hspace{2in} &  & \vdots\\
\tr(\sym^{k-n+2} A) - \tr(A)\tr(\sym^{k-n+1} A) + \tr(\ext^2 A)\tr(\sym^{k-n} A) +\cdots  
+(-1)^n\tr(\ext^n A)\tr(\sym^{k-2n+2} A) & = & 0
\end{array}
\right. 
\]
}
Multiply them by $(-1)^0\sigma_1(A), (-1)^1\sigma_2(A),\ldots,(-1)^{n-2}\sigma_{n-1}(A)$ respectively 
{\footnotesize
\[
\left\{
\begin{array}{rcc}
\sigma_1(A)[\tr(\sym^{k} A) - \tr(A)\tr(\sym^{k-1} A)  +\cdots   
+(-1)^n\tr(\ext^n A)\tr(\sym^{k-n} A)] & = & 0\\
-\sigma_2(A)[\tr(\sym^{k-1} A) - \tr(A)\tr(\sym^{k-2} A)  +\cdots  
+(-1)^n\tr(\ext^n A)\tr(\sym^{k-n-1} A)] & = & 0\\
\vdots\hspace{2in} &  & \vdots\\
(-1)^{n-2}\sigma_{n-1}(A)[\tr(\sym^{k-n+2} A) - \tr(A)\tr(\sym^{k-n+1} A) 
+\cdots 
+(-1)^n\tr(\ext^n A)\tr(\sym^{k-2n+2} A)] & = & 0
\end{array}
\right. 
\]
}
and add them all up to get
{\footnotesize
\begin{eqnarray}
0
&=&[\sigma_1(A)\tr(\sym^{k} A) - \sigma_2(A)\tr(\sym^{k-1} A)
+\ldots+(-1)^{n-2}\sigma_{n-1}(A)\tr(\sym^{k-n+2}A)]\nonumber\\
&&-[\sigma_1(A)\tr(\sym^{k-1} A) - \sigma_2(A)\tr(\sym^{k-2} A)
+\ldots+(-1)^{n-2}\sigma_{n-1}(A)\tr(\sym^{k-n+1}A)]\tr(A)\nonumber\\
&&+[\sigma_1(A)\tr(\sym^{k-2} A) - \sigma_2(A)\tr(\sym^{k-3} A)
+\ldots+(-1)^{n-4}\sigma_{n-1}(A)\tr(\sym^{k-n}A)]\tr(\ext^2A)\label{eq:intermediate}\\
&&\vdots \nonumber\\
&&+(-1)^n[\sigma_1(A)\tr(\sym^{k-n} A) - \sigma_2(A)\tr(\sym^{k-n-1} A)
+\ldots+(-1)^{n-2}\sigma_{n-1}(A)\tr(\sym^{k-2n+2}A)]\tr(\ext^nA).\nonumber
\end{eqnarray}
}
In order to simplify the notation, let us define
\[P^{k+1} (A) := \sigma_1(A)\tr(\sym^{k} A) - \sigma_2(A)\tr(\sym^{k-1} A)
+\ldots+(-1)^{n-2}\sigma_{n-1}(A)\tr(\sym^{k-n+2}A),\]
so that equation \rf{eq:intermediate} becomes
\[P^{k+1} (A)-\tr(A)P^{k} (A) +\tr(\ext^2 A)P^{k-1}(A)+\cdots
+(-1)^n\tr(\ext^n A)P^{k-n+1}(A)= 0\]
which is the same equation as the one fulfilled by $\sigma_1(A^k)$ due to the Cayley-Hamilton
theorem
and the linearity of $\sigma_1$
\[\sigma_1(A^{k+1}) -  \tr( A) \sigma_1(A^k)  +  \tr( \ext^2 A) \sigma_1(A^{k-1})  
%- \tr(\ext^3 A)\sigma_1(A^{k-2}) 
+\cdots + (-1)^n\tr(\ext^n A)\sigma_1(A^{k-n+1})
=0,\]
where $k\geq n-1$.
Since they are recursive formulas, all we have to do is check that the first $n-1$ terms coincide,
i.e.
we have to check that 
\begin{eqnarray*}
\sigma_1(A^0)&=& P^0 (A),\\ 
\sigma_1(A^1)&=& P^1 (A),\\ 
\vdots && \vdots \\
\sigma_1(A^{n-1})&=& P^{n-1} (A). 
\end{eqnarray*}
To begin with, recall that $\tr(\sym^{-s} A)=0$ for $s\geq 1$, so that $\sigma_1=P^1$.
Next, the identities for $\sigma_2(A),\ldots,\sigma_{n-1}(A)$ are exactly the recursive definition
of such functions.
Hence
\[\sigma_1(A^{k+1}) = \sigma_1(A)\tr(\sym^{k} A) - \sigma_2(A)\tr(\sym^{k-1} A)
+\ldots+(-1)^{n-2}\sigma_{n-1}(A)\tr(\sym^{k-n+2}A)\]
for arbitrary $k \in\mathbb{Z}$. 
Thus, we have proved the following.

\begin{theo}\label{thm:general-case} Let $K$ be a field, $A\in M_{n\times n}(K)$, $\sigma_1$ an arbitrary
linear funtional on $M_{n\times n}$ which vanishes on ${\rm Id}_{n\times n}$, and polynomial functions 
$\sigma_j$, $2\leq j\leq n-1$ defined recursively in {\em \rf{eq:sigmas}}. Then 
\[\sigma_1(A^{k+1}) = \sigma_1(A)\tr(\sym^{k} A) - \sigma_2(A)\tr(\sym^{k-1} A)
+\ldots+(-1)^{n-2}\sigma_{n-1}(A)\tr(\sym^{k-n+2}A)\]
for arbitrary $k \in\mathbb{Z}$. 
% with $\sigma_j$, $1\leq j\leq n-1$, as defined in {\em \rf{eq:sigmas}}.
Furthermore, if $\sigma_1(A)\not=0$
\[\tr(\sym^{k} A)={1\over \sigma_1(A)}[\sigma_1(A^{k+1})  + \sigma_2(A)\tr(\sym^{k-1} A)
+\ldots+(-1)^{n-1}\sigma_{n-1}(A)\tr(\sym^{k-n+2}A)].\]
\end{theo}
\qd

This means that any $\tr(\sym^kA)$ can be written in terms of values $\sigma_i(A^j)$, where
$1\leq j\leq k+1$, $1\leq i\leq n-1$.

First, we recover the result in \cite{Cisneros}.

\begin{corol} For a $2\times 2$ matrix $A$
\[\sigma_1(A^{k+1}) = \sigma_1(A)\tr(\sym^{k} A) \]
for arbitrary $k \in\mathbb{Z}$. 
If $\sigma_1(A)\not=0$
\[\tr(\sym^{k} A)={\sigma_1(A^{k+1})\over \sigma_1(A)}.\]
\end{corol}
\qd

One can deduce related formulas for $n=3,4$.

\begin{corol} For a $3\times 3$ matrix $A$, if $\sigma_1(A)\not=0$,
\begin{eqnarray*}
\tr(\sym^k A) 
&=& {1\over \sigma_1(A)}
\left[\sum_{i=0}^{k-2}\left({\sigma_2(A)\over\sigma_1(A)}\right)^i
\sigma_1(A^{k+1-i})\right] + \left({\sigma_2(A)\over\sigma_1(A)}\right)^{k-1}
{\sigma_1(A^2)+\sigma_2(A)\over\sigma_1(A)} .
\end{eqnarray*}
\end{corol}
\qd

\begin{corol} For a $4\times 4$ matrix $A$, if $\sigma_1(A)\not=0$,
{\footnotesize
\begin{eqnarray*}
\tr(\sym^k A)
   &=& {1\over \sigma_1(A)} 
   \left[\sum_{j=0}^{[{k-3\over2}]} \left({\sigma_3(A)\over
\sigma_1(A)}\right)^j(-1)^j\sum_{i=0}^{k-(2j+3)}
   {i+j\choose j}\left({\sigma_2(A)\over \sigma_1(A)}\right)^i \sigma_1(A^{k+1-i-2j})    \right]\\
   && +\left( {\sigma_1(A^3)-\sigma_3(A)\over\sigma_1(A)} + {\sigma_2(A)\over\sigma_1(A)}{
\sigma_1(A^2)+\sigma_2(A)\over \sigma_1(A)} \right) 
\sum_{j=0}^{[{k-3\over2}]}(-1)^j{k-2-j\choose j}\left({\sigma_3(A)\over
\sigma_1(A)}\right)^j
      \left({\sigma_2(A)\over \sigma_1(A)}\right)^{k-2-2j}\\
   && -{\sigma_1(A^2)+\sigma_2(A)\over\sigma_1(A)} \left({\sigma_3(A)\over
\sigma_1(A)}\right)\sum_{j=0}^{[{k-3\over2}]}(-1)^j{k-3-j\choose j}\left({\sigma_3(A)\over
\sigma_1(A)}\right)^j
      \left({\sigma_2(A)\over \sigma_1(A)}\right)^{k-3-2j}.
\end{eqnarray*}
}
\end{corol}
\qd

\section{Examples}\label{sec:examples}

\subsection{Case $n=3$}

Setting $\sigma_1$ equal to the various linear functionals of $3\times 3$ matrices 
vanishing on the identity matrix, we have the accompanying functions $\sigma_2$ as follows:
\[
\begin{array}{|c|c|}\hline
\sigma_1(A) & \sigma_2(A)\\\hline\hline
a_{12} & a_{12}a_{33}-a_{13}a_{32} \\\hline
a_{13} & a_{13}a_{22}-a_{12}a_{23}\\\hline
a_{23} & a_{11}a_{23}-a_{21}a_{13}\\\hline
a_{21} & a_{21}a_{33}-a_{23}a_{31}\\\hline
a_{31} & a_{31}a_{22}-a_{32}a_{21}\\\hline
a_{32} & a_{32}a_{11}-a_{31}a_{12}\\\hline
a_{11}-a_{22} & a_{11}a_{33}-a_{22}a_{33}-a_{13}a_{31}+a_{23}a_{32}\\\hline
a_{11}-a_{33} & a_{11}a_{22}-a_{22}a_{33}-a_{12}a_{21}+a_{23}a_{32}\\\hline
\end{array} 
\]

\subsection{Case $n=4$}

Setting $\sigma_1$ equal to the various linear functionals of $3\times 3$ matrices 
vanishing on the identity matrix, we have the accompanying functions $\sigma_2$ and $\sigma_3$ as
follows:
\[
\begin{array}{|c|c|c|}\hline
\sigma_1(A) & \sigma_2(A) & \sigma_3(A)\\\hline\hline
a_{12} & a_{12}a_{33}+a_{12}a_{44}-a_{13}a_{32}-a_{14}a_{42} &
a_{12}a_{33}a_{44}-a_{12}a_{34}a_{43}-a_{32}a_{13}a_{44}\\ &&
+a_{32}a_{14}a_{43}+a_{42}a_{13}a_{34}-a_{42 }a_{14}a_{33}\\\hline
a_{13} & a_{13}a_{22}+a_{13}a_{44}-a_{12}a_{23}-a_{14}a_{43} &
-a_{12}a_{23}a_{44}+a_{12}a_{24}a_{43}+a_{22}a_{13}a_{44}\\
&&-a_{22}a_{14}a_{43}-a_{42}a_{13}a_{24}+a_{4 2}a_{14}a_{23}\\\hline
a_{14} & a_{14}a_{22}+a_{14}a_{33}-a_{12}a_{24}-a_{13}a_{34} &
-a_{22}a_{13}a_{34}+a_{22}a_{14}a_{33}+a_{32}a_{13}a_{24}\\
&&-a_{32}a_{14}a_{23}+a_{12}a_{23}a_{34}-a_{1 2}a_{24}a_{33}\\\hline
a_{23} & a_{23}a_{11}+a_{23}a_{44}-a_{21}a_{13}-a_{24}a_{43} &
a_{43}a_{21}a_{14}+a_{13}a_{24}a_{41}+a_{23}a_{11}a_{44}\\
&&-a_{23}a_{14}a_{41}-a_{21}a_{13}a_{44}-a_{24 }a_{43}a_{11}\\\hline
a_{24} & a_{24}a_{11}+a_{24}a_{33}-a_{21}a_{14}-a_{23}a_{34} &
-a_{24}a_{13}a_{31}+a_{14}a_{23}a_{31}+a_{34}a_{21}a_{13}\\
&&+a_{24}a_{11}a_{33}-a_{21}a_{14}a_{33}-a_{2 3}a_{34}a_{11}\\\hline
a_{34} & a_{34}a_{11}+a_{34}a_{22}-a_{31}a_{14}-a_{32}a_{24} &
-a_{34}a_{12}a_{21}+a_{14}a_{32}a_{21}+a_{34}a_{11}a_{22}\\
&&+a_{24}a_{31}a_{12}-a_{31}a_{14}a_{22}-a_{3 2}a_{24}a_{11}\\\hline
a_{21} & a_{21}a_{33}+a_{21}a_{44}-a_{23}a_{31}-a_{24}a_{41} &
a_{31}a_{24}a_{43}-a_{21}a_{34}a_{43}+a_{41}a_{23}a_{34}\\
&&+a_{21}a_{33}a_{44}-a_{23}a_{31}a_{44}-a_{24 }a_{41}a_{33}\\\hline
a_{31} & a_{31}a_{22}+a_{31}a_{44}-a_{32}a_{21}-a_{34}a_{41} &
a_{41}a_{32}a_{24}+a_{21}a_{34}a_{42}+a_{31}a_{22}a_{44}\\
&&-a_{31}a_{24}a_{42}-a_{32}a_{21}a_{44}-a_{34 }a_{41}a_{22}\\\hline
a_{41} & a_{41}a_{22}+a_{41}a_{33}-a_{42}a_{21}-a_{43}a_{31} &
a_{41}a_{22}a_{33}+a_{21}a_{43}a_{32}+a_{31}a_{42}a_{23}\\
&&-a_{41}a_{23}a_{32}-a_{42}a_{21}a_{33}-a_{43 }a_{31}a_{22}\\\hline
a_{32} & a_{32}a_{11}+a_{32}a_{44}-a_{31}a_{12}-a_{34}a_{42} &
a_{12}a_{34}a_{41}+a_{42}a_{31}a_{14}+a_{32}a_{11}a_{44}\\
&&-a_{32}a_{14}a_{41}-a_{31}a_{12}a_{44}-a_{34 }a_{42}a_{11}\\\hline
a_{42} & a_{42}a_{11}+a_{42}a_{33}-a_{41}a_{12}-a_{43}a_{32} &
a_{12}a_{43}a_{31}+a_{32}a_{41}a_{13}+a_{42}a_{11}a_{33}\\
&&-a_{42}a_{13}a_{31}-a_{41}a_{12}a_{33}-a_{43 }a_{32}a_{11}\\\hline
a_{43} & a_{43}a_{11}+a_{43}a_{22}-a_{41}a_{13}-a_{42}a_{23} &
a_{43}a_{11}a_{22}-a_{43}a_{12}a_{21}+a_{13}a_{42}a_{21}\\
&&+a_{23}a_{41}a_{12}-a_{41}a_{13}a_{22}-a_{42 }a_{23}a_{11}\\\hline
a_{11}-a_{22} &
a_{11}a_{33}-a_{22}a_{33}+a_{11}a_{44}-a_{22}a_{44} &
-a_{14}a_{41}a_{33}-a_{13}a_{31}a_{44}-a_{34}a_{43}a_{11} \\
& -a_{13}a_{31}-a_{14}a_{41}+a_{23}a_{32}+a_{24}a_{42 }
&+a_{33}a_{44}a_{11}+a_{31}a_{14}a_{43}+a_{41}a_{13}a_{34}\\
&&-a_{42}a_{23}a_{34}+a_{42}a_{24}a_{33}-a_{22}a_{33}a_{44}\\
&&+a_{22}a_{34}a_{43}+a_{32}a_{ 23}a_{44}-a_{32}a_{24}a_{43}\\\hline
a_{11}-a_{33} &
a_{11}a_{22}-a_{22}a_{33}+a_{11}a_{44}-a_{33}a_{44} & 
a_{21}a_{14}a_{42}-a_{24}a_{42}a_{11}-a_{14}a_{41}a_{22}\\
&-a_{12}a_{21}-a_{14}a_{41}+a_{23}a_{32}+a_{34}a_{43}
&+a_{11}a_{44}a_{22}+a_{41}a_{12}a_{24}-a_{42}a_{23}a_{34}\\
&&+a_{42}a_{24}a_{33}-a_{22}a_{33}a_{44}+a_{22}a_{34}a_{43}\\
&&+a_{32}a_{23}a_{44}-a_{32}a_{2 4}a_{43}-a_{12}a_{21}a_{44}\\\hline
a_{11}-a_{44} &
a_{11}a_{22}-a_{22}a_{44}+a_{11}a_{33}-a_{33}a_{44}
&a_{31}a_{12}a_{23}+a_{21}a_{13}a_{32}-a_{23}a_{32}a_{11}\\
&-a_{12}a_{21}-a_{13}a_{31}+a_{24}a_{42}+a_{34}a_{43}
&-a_{13}a_{31}a_{22}+a_{11}a_{33}a_{22}-a_{42}a_{23}a_{34}\\
&&+a_{42}a_{24}a_{33}-a_{22}a_{33}a_{44}+a_{22}a_{34}a_{43}\\
&&+a_{32}a_{23}a_{44}-a_{32}a_{2 4}a_{43}-a_{12}a_{21}a_{33}\\
\hline
\end{array}
\]

{\small 
\renewcommand{\baselinestretch}{0.5}
\newcommand{\bi}{\vspace{-.05in}\bibitem} }

\enddocument